\documentclass[preprint,12pt]{elsarticle}

\usepackage{amssymb}
\usepackage{amsthm}
\usepackage{enumerate}

\newtheorem{defn}{Definition}
\newtheorem{prop}{Proposition}
\newtheorem{exmp}{Example}
\newtheorem{rem}{Remark}
\newproof{pf}{Proof}
\newproof{ack}{Acknowledgements}

\begin{document}

\begin{frontmatter}


\cortext[cor1]{Corresponding author}

\title{Intuitionistic fuzzy parametrized soft set theory and its decision making}


\author[rvt1]{Irfan Deli\corref{cor1}}
\ead{irfandeli@kilis.edu.tr}
\author[rvt2]{{Naim \c{C}a\u{g}man}}
\ead{ncagman@gop.edu.tr}
\address[rvt1]{Department of
Mathematics, Faculty of Arts and Sciences, 7 Aral{\i}k University ,
79000 Kilis, Turkey}
\address[rvt2]{Department of Mathematics, Faculty of Arts and Sciences,
Gaziosmanpa\c{s}a University 60250 Tokat, Turkey}

\begin{abstract}
In this work, we first define intuitionistic fuzzy parametrized soft
sets (intuitionistic FP-soft sets) and study some of their
properties.  We then introduce an adjustable approaches to
intuitionistic FP-soft sets based decision making. We also give an
example which shows that they can be successfully applied to
problems that contain uncertainties.
\end{abstract}

\end{frontmatter}


\section{Introduction}
Many fields deal with the uncertain data which may not be
successfully modeled by the classical mathematics, probability
theory, fuzzy sets \cite{zad-65}, rough sets \cite{paw-82}, and
other mathematical tools. In 1999, Molodtsov \cite{mol-99} proposed
a completely new approach so-called {\em soft set theory} that is
more universal for modeling vagueness and uncertainty.

After definitions of the operations of soft sets
\cite{ali-09,cag-09a,maj-03,aslý-11a}, the properties and
applications on the soft set theory have been studied increasingly
(e.g. \cite{ali-09,cag-11d,cag-09b}). The algebraic structure of
soft set theory has also been studied increasingly (e.g.
\cite{acar-10,akt-07,aslý-11b,fen-08,jun-10a,jun-08a,jun-08b,zhan-10,par-08}).
In recent years, many interesting applications of soft set theory
have been expanded by embedding the ideas of fuzzy sets (e.g. \cite{
ali-11,ayg-09,cag-09d, fen-08b,zhan-10,roy-07,maj-01a,maju-10}),
rough sets (e.g. \cite{ali-11, fen-11, fen-10}) and intuitionistic
fuzzy sets (e.g. \cite{maj-01c, maj-01b,muk-08})

\c{C}a\u{g}man \emph{et al.}\cite{cag-09c} defined FP-soft sets and
constructed an FP-soft set decision making method. In this paper, we
first define intuitionistic fuzzy parametrized soft sets
(intuitionistic FP-soft sets), and study their operations and
properties. We then introduce a decision making method based on
intuitionistic FP-soft sets. This method is more practical and can
be successfully applied to many problems that contain uncertainties.
\section{Preliminary}\label{ss}
In this section, we present the basic definitions of soft set theory
\cite{mol-99}, fuzzy set theory \cite{zad-65}, intuitionistic fuzzy
set theory \cite{ata-86} and FP-soft set theory \cite{cag-09c} that
are useful for subsequent discussions.
\begin{defn}
Let $U$ be an initial universe,  $P(U)$ be the power set of $U$, $E$
is a set of parameters and $A \subseteq E$. Then, a soft set $F_A$
over $U$ is defined as follows:
$$
F_A= \{(x, f_A(x)): x\in E\}
$$
where $f_A: E\to P(U)$ such that $f_A(x)=\emptyset$ if $x\notin A$.

Here, $f_A$ is called approximate function of the soft set $F_A$,
and the value $f_A(x)$ is a set called \emph{$x$-element} of the
soft set for all  $x \in E$. It is worth noting that the sets
$f_A(x)$ may be arbitrary.
\end{defn}
\begin{defn}
Let $U$ be a universe. Then a fuzzy set $X$ over $U$ is a function
defined as follows:
$$
X=\{(\mu_X(u)/u): u\in U\}
$$
where $\mu_X:U\rightarrow [0.1]$

Here, $\mu_X$ called membership function of $X$, and the value
$\mu_X(u)$ is called the grade of membership of $u\in U $. The value
represents the degree of u belonging to the fuzzy set $X$.
\end{defn}
\begin{defn}
Let $E$ be a universe. An intuitionistic fuzzy set $A$ on $E$ can be
defined as follows:
$$
A =\{<x, \mu_{A}(x),\gamma_{A}(x)>: \,\, x \in E \}
$$
where, $\mu_{A}: E \rightarrow [ 0, 1 ]$ and $\gamma_{A} : E
\rightarrow [ 0, 1 ]$ such that $0 \leq \mu_{A} (x) + \gamma_A(x)
\leq 1$ for any $ x\in E$.

Here, $\mu_{A}(x)$ and $\gamma_{A}(x)$ is the degree of membership
and degree of non-membership of the element $x$, respectively.
\end{defn}

If $A$ and $B$ are two intuitionistic fuzzy sets on $E$, then
\begin{enumerate}[i.]
\item $A \subset B$ if and only if
$\mu_{A}(x)\leq \mu_{B}(x)$ and $\gamma_{A}(x) \geq \gamma_{B}(x)$
for $\forall x \in E$

\item $A=B$ if and only if  $\mu_{A}(x)=
\mu_{B}(x)$ and $\gamma_{A}(x) = \gamma_{B}(x)$ $\forall x \in E$

\item $A ^c= \{<x,  \gamma_{A}(x), \mu_{A}(x)>:\,\, x \in E \}$

\item $A\cup B = \{<x, max(\mu_{A}(x),\mu_{B}(x)), min(  \gamma_{A}(x),
\gamma_{B}(x)>:\,\, x \in E \} $,

\item $A\cap B = \{<x, min(\mu_{A}(x),\mu_{B}(x)), max(  \gamma_{A}(x),
\gamma_{B}(x)>:\,\, x \in E \} $,

\item $A + B = \{<x,
\mu_{X}(x)+ \mu_{Y}(x)- \mu_{X}(x)\mu_{Y}(x),
\gamma_{X}(x)\gamma_{Y}(x)>:\,\, x \in E \}$,

\item $A\cdot B
= \{<x, \mu_{A}(x)\mu_{B}(x), \gamma_{A}(x) + \gamma_{B}(x)-
\gamma_{A}(x)\gamma_{B}(x)>:\,\, x \in E \}$.
\end{enumerate}
\begin{defn}
Let $U$ be an initial universe,  $P(U)$ be the power set of $U$, $E$
be a set of all parameters and $X$ be a fuzzy set over $E$. Then a
FP-soft set $(f_X, E)$ on the universe $U$ is defined as follows:
$$
\label{FP-soft set} (f_X, E)= \{(\mu_X(x)/x, f_X(x)): x\in E\}
$$
where $\mu_X:E\rightarrow [0.1]$ and $f_X: E\to P(U)$ such that
$f_X(x)=\emptyset$ if $\mu_X(x)=0$.

Here  $f_X$ called approximate function and $\mu_X$ called
membership function of FP-soft sets.
\end{defn}
\section{Intuitionistic FP-soft sets}
In this section, we define intuitionistic fuzzy parametrized soft
sets (intuitionistic FP-soft sets) and their operations.
\begin{defn}
 Let $U$ be an initial universe,  $P(U)$ be the power set of~ $U$,
 $E$ be a set of all parameters and $K$ be an
intuitionistic fuzzy set over $E$. An intuitionistic FP-soft sets
$\amalg_K$ over U is defined as follows:
$$
\begin{array}{rl}
\label{IFP-soft set} \amalg_K = &\bigg\{(<x,\alpha_K(x),\beta_K(x)
>, f_K(x)): x\in E\bigg\}
\end{array}
$$
where $\alpha_K:E\rightarrow [0.1]$, $\beta_K:E\rightarrow [0.1]$
and $f_K: E\to P(U)$ with the property $f_K(x)=\emptyset$ if
$\alpha_K(x)=0$ and $\beta_K(x)=1$.

Here, the function $\alpha_K$ and $\beta_K$ called membership
function and non-membership of intuitionistic FP-soft set,
respectively. The value $\alpha_K(x)$ and $\beta_K(x)$ is the degree
of importance and unimportant of the parameter $x$.
\end{defn}


Obviously, each ordinary FP-soft set can be written as
$$
\label{IFP-soft set} \amalg_K = \bigg\{(<x,
\alpha_K(x),1-\alpha_K(x)
>, f_K(x)): x\in E\bigg\}
$$
Note that the sets of all intuitionistic FP-soft sets over $U$ will
be denoted by $IFPS(U) $.
\begin{defn}
Let $\amalg_K\in IFPS(U) $. If $\alpha_K(x)=0$ and $\beta_K(x)=1$
for all $x\in E$, then $\amalg_K$ is called a empty intuitionistic
FP-soft sets, denoted by $\amalg_\Phi$.
\end{defn}
\begin{defn}
Let $\amalg_K\in IFPS(U) $. If  $\alpha_K(x)=1$, $\beta_K(x)=0$ and
$f_K(x)=U$ for all $x\in E $, then $\amalg_K$ is called universal
intuitionistic FP-soft set, denoted by $\amalg_{\tilde{E}}$.
\end{defn}
\begin{exmp}\label{ex-soft}
Assume that $U=\{u_1, u_2, u_3, u_4, u_5\}$ is a universal set and
$E=\{x_1,x_2, x_3,x_4\}$ is a set of parameters. If
$$
K=\{<x_2,0.2,0.5 >, <x_3,0.5,0.5 >, <x_4,0.6,0.3 >\}
$$
and
$$
f_K(x_2)= \{u_2, u_4\}, f_K(x_3)= \emptyset,  f_K(x_4)= U
$$
then an intuitionistic FP-soft set $\amalg_K$ is written by
$$
\amalg_K=\{(<x_2,0.2,0.5 >, \{u_2, u_4\}), (<x_3,0.5,0.5
>,\emptyset), (<x_4,0.6,0.3 >,U)\}
$$

If $L=\{<x_1,0,1>, <x_2,0,1>, <x_3,0,1>, <x_4,0,1>\}$, then the
intuitionistic FP-soft set $\amalg_L$ is an empty intuitionistic FP
soft set.

If $M= \{<x_1,1,0>, <x_2,1,0>,<x_3,1,0>, <x_4,1,0>\}$ and $
f_M(x_1)=U$, $ f_M(x_2)= U$,$f_M(x_3)=U$ and $ f_M(x_4)= U$, then
the intuitionistic FP-soft set $\amalg_M$ is a universal
intuitionistic FP-soft set.
\end{exmp}
\begin{defn}
$\amalg_K, \amalg_L \in IFPS(U) $. Then $\amalg_K$ is a
intuitionistic FP-soft subset of $\amalg_L$, denoted by $\amalg_K
\widetilde{\subseteq} \amalg_L$, if and only if $\alpha_K(x) \leq
\alpha_L(x)$, $\beta_K(x) \geq \beta_L(x) \textrm{ and }
f_K(x)\subseteq f_L(x)$ for all $x \in E$.
\end{defn}
\begin{rem}
$\amalg_K \widetilde{\subseteq} \amalg_L$ does not imply that every
element of $\amalg_K$ is an element of $\amalg_L$ as in the
definition of classical subset. For example, assume that $U=\{u_{1},
u_{2}, u_{3}, u_{4}\}$ is a universal set of objects and $E=\{x_{1},
x_{2}, x_{3}\}$ is a set of all parameters. If $K= \{<x_1,0.4,0.6
>\}$ and $L=\{<x_1,0.5,0.5>,$ $<x_3,0.4,0.5 >\}$, and
$\amalg_{K}=\{(<x_1,0.4,0.6 >, \{u_2, u_4\})\}$,
$\amalg_{L}=\{(<x_1,0.5,0.5 >, \{u_2, u_3, u_4\}),(<x_3,0.4,0.5
>, \{u_1, u_5\})\}$, then for all $x \in E $, $\alpha_K(x) \leq
\alpha_L(x)$, $\beta_K(x) \geq \beta_L(x) \textrm{ and }
\amalg_K(x)\subseteq \amalg_L(x)$ is valid. Hence $\amalg_K
\widetilde{\subseteq} \amalg_L $. It is clear that $(<x_1,0.4,0.6
>,\{u_2,u_4\}) \in \amalg_K $ but $(<x_1,0.4,0.6>,$ $\{u_2,u_4\})\notin \amalg_L$.
\end{rem}
\begin{prop} Let $\amalg_K,\amalg_L\in IFPS(U) $. Then
\begin{enumerate}[i.]
\item $\amalg_K \widetilde{\subseteq} \amalg_{\tilde E}$
\item $\amalg_\Phi \widetilde{\subseteq} \amalg_K $
\item $\amalg_K\widetilde{\subseteq} \amalg_K $
\end{enumerate}
\end{prop}
\begin{defn}
$\amalg_K, \amalg_L \in IFPS(U) $. Then $\amalg_K$ and $\amalg_L$
are intuitionistic FP-soft-equal, written by $\amalg_K = \amalg_L$,
if and only if $\alpha_K(x) = \alpha_L(x)$ , $\beta_K(x)=
\beta_L(x)$  and $f_K(x)= f_L(x)$ for all $x \in E$.
\end{defn}

\begin{prop} Let $\amalg_K, \amalg_L , \amalg_M\in IFPS(U) $. Then
\begin{enumerate}[i.]
\item $\amalg_K=\amalg_L$ and $\amalg_L = \amalg_M  \Leftrightarrow \amalg_K = \amalg_M$
\item $\amalg_K \widetilde{\subseteq} \amalg_L$ and $\amalg_L \widetilde{\subseteq} \amalg_K  \Leftrightarrow \amalg_K = \amalg_L$
\item $\amalg_K \widetilde{\subseteq} \amalg_L$ and $\amalg_L
\widetilde{\subseteq} \amalg_M  \Rightarrow \amalg_K
\widetilde{\subseteq} \amalg_M$
\end{enumerate}
\end{prop}
\begin{defn}
$\amalg_K \in IFPS(U) $. Then complement of $\amalg_K$, denoted by
$\amalg_K^c$, is a intuitionistic FP-soft set defined by
$$
\label{IFP-soft set} \amalg_{K}^c = \bigg\{(<x,
\beta_K(x),\alpha_K(x)>, f_{K^c}(x)): x\in K\bigg\}
$$
where $f_{K^c}(x)= U\setminus f_K(x)$.
\end{defn}
\begin{prop}
Let $\amalg_K\in IFPS(U) $. Then
\begin{enumerate}[i.]
\item $(\amalg_K^c)^{c}= \amalg_K$
\item $ \amalg_\Phi^c = \amalg_{\tilde E}$
\item $  \amalg_{\tilde E}^c= \amalg_\Phi$
\end{enumerate}
\end{prop}
\begin{defn}\label{union}
$\amalg_K, \amalg_L \in IFPS(U) $. Then union of $\amalg_K$ and
$\amalg_L$, denoted by $ \amalg_K \widetilde{\cup} \amalg_L$,  is
defined by
$$
\amalg_K{\widetilde{\cup}\amalg_L}=\bigg\{(<x,
max(\alpha_{K}(x),\alpha_{L}(x)), min( \beta_{K}(x),
\beta_{L}(x))>,f_{K \widetilde{\cup} L}(x)): x \in E \bigg\}
$$
where $f_{K \widetilde{\cup} L}(x)=f_K(x)\cup f_L(x)$.
\end{defn}
\begin{prop}
Let $\amalg_K,\amalg_L,\amalg_M\in IFPS(U) $. Then
\begin{enumerate}[i.]
\item  $ \amalg_K \widetilde{\cup} \amalg_K = \amalg_K $
\item $ \amalg_K \widetilde{\cup} \amalg_\Phi = \amalg_K$
\item  $ \amalg_K \widetilde{\cup} \amalg_{\tilde E} = \amalg_{\tilde E} $
\item $ \amalg_K \widetilde{\cup} \amalg_L= \amalg_L\widetilde{\cup} \amalg_K $
\item $(\amalg_K \widetilde{\cup} \amalg_L)\widetilde{\cup} \amalg_M= \amalg_K \widetilde{\cup} (\amalg_L\widetilde{\cup} \amalg_M) $
\end{enumerate}
\end{prop}

\begin{defn}\label{intersection}
$\amalg_K, \amalg_L \in IFPS(U) $. Then intersection of $\amalg_K$
and $\amalg_L$, denoted by $ \amalg_K \widetilde{\cap} \amalg_L$, is
a intuitionistic FP-soft sets defined by
$$
\amalg_K{\widetilde{\cap}\amalg_L}=\bigg\{<x,
min(\alpha_{K}(x),\alpha_{L}(x)), max(  \beta_{K}(x),
\beta_{L}(x)>,f_{K\widetilde{\cap} L}(x)):  x \in E \bigg\}
$$
where $f_{K\widetilde{\cap} L}(x)= f_K(x)\cap f_L(x)$.
\end{defn}
\begin{prop}
Let $\amalg_K,\amalg_L,\amalg_M\in IFPS(U) $. Then
\begin{enumerate}[i.]
\item $ \amalg_K \widetilde{\cap} \amalg_K = \amalg_K
$
\item $ \amalg_K \widetilde{\cap} \amalg_\Phi = \amalg_\Phi$
\item $ \amalg_K \widetilde{\cap} \amalg_{\tilde E} = \amalg_K  $
\item $ \amalg_K \widetilde{\cap} \amalg_L= \amalg_L\widetilde{\cap} \amalg_K $
\item $(\amalg_K \widetilde{\cap} \amalg_L)\widetilde{\cap} \amalg_M= \amalg_K \widetilde{\cap} (\amalg_L\widetilde{\cap} \amalg_M) $
\end{enumerate}
\end{prop}

\begin{rem}
Let $\amalg_K\in IFPS(U) $. If $\amalg_K\neq \amalg_\Phi $ or
$\amalg_K\neq \amalg_{\tilde E}$, then $\amalg_K \widetilde{\cup}
\amalg_K^{c}\neq \amalg_{\tilde E}$ and $\amalg_K \widetilde{\cap}
\amalg_K^{c}\neq \amalg_{\Phi}$.
\end{rem}
\begin{prop}
Let $\amalg_K,\amalg_L,\amalg_M\in IFPS(U) $. Then
\begin{enumerate}[i.]
\item $ \amalg_K \widetilde{\cup} (\amalg_L \widetilde{\cap} \amalg_M)= (\amalg_K \widetilde{\cup} \amalg_L) \widetilde{\cap} (\amalg_K \widetilde{\cup} \amalg_M) $
\item $ \amalg_K \widetilde{\cap} (\amalg_L \widetilde{\cup} \amalg_M)= (\amalg_K\widetilde{\cap} \amalg_L)\widetilde{\cup} (\amalg_K \widetilde{\cap} \amalg_M) $
\end{enumerate}
\item
\end{prop}

\begin{prop}
Let $\amalg_K,\amalg_L IFPS(U) $. Then following DeMorgan's types of
results are true.
\begin{enumerate}[{i.}]
\item   $(\amalg_K \widetilde{\cup} \amalg_L )^{c} =\amalg_K ^{c} \widetilde{\cap}\amalg_L ^{c}$
\item   $(\amalg_K \widetilde{\cap}  \amalg_L )^{c} =\amalg_K ^{c} \widetilde{\cup} \amalg_L ^{c}$
\end{enumerate}
\end{prop}
\begin{defn}\label{or-sum}
$\amalg_K, \amalg_L \in IFPS(U) $. Then OR-sum of $\amalg_K$ and
$\amalg_L$, denoted by $ \amalg_K\vee^{+} \amalg_L$,  is defined by
$$
\begin{array}{l}
{\amalg_K \vee^{+} \amalg_L}(x)= \\
\bigg\{(<x,\alpha_K(x)+\alpha_L(x)-\alpha_K(x)\alpha_L(x),
\beta_K(x)\beta_L(x)>,f_{K \widetilde{\cup} L}(x)): x \in E \bigg\}
\end{array}
$$
where $f_{K \widetilde{\cup} L}(x)= f_K(x)\cup f_L(x)$.
\end{defn}
\begin{defn}\label{and-sum}
$\amalg_K, \amalg_L \in IFPS(U) $. Then AND-sum of $\amalg_K$ and
$\amalg_L$, denoted by $ \amalg_K\wedge^{+}  \amalg_L$,  is defined
by
$$
\begin{array}{l}
{\amalg_K \wedge^{+}
\amalg_L}(x)=\\\bigg\{(<x,\alpha_K(x)+\alpha_L(x)-\alpha_K(x)\alpha_L(x),
\beta_K(x)\beta_L(x)
>,f_{K \widetilde{\cap} L}(x)): x \in E \bigg\}
\end{array}
$$
where $f_{K \widetilde{\cap} L}(x)= f_K(x)\cap f_L(x)$.
\begin{prop}
Let $\amalg_K,\amalg_L,\amalg_M\in IFPS(U) $. Then
\begin{enumerate}[i.]
\item $ \amalg_K \vee^{+}  \amalg_\Phi = \amalg_K$
\item $ \amalg_K \vee^{+} \amalg_{\tilde E} = \amalg_{\tilde E} $
\item $ \amalg_K\vee^{+}  \amalg_L= \amalg_L\vee^{+} \amalg_K $
\item $ \amalg_K \wedge^{+} \amalg_L= \amalg_L\wedge^{+} \amalg_K $
\item $(\amalg_K \vee^{+}  \amalg_L)\vee^{+} \amalg_M= \amalg_K \vee^{+}  (\amalg_L\vee^{+}  \amalg_M) $
\item $(\amalg_K \wedge^{+} \amalg_L)\wedge^{+} \amalg_M= \amalg_K \wedge^{+} (\amalg_L\wedge^{+} \amalg_M) $
\end{enumerate}
\end{prop}
\end{defn}
\begin{defn}\label{or-pro}
$\amalg_K, \amalg_L \in IFPS(U) $. Then OR-product $\amalg_K$ and
$\amalg_L$, denoted by $ \amalg_K \vee^{\times} \amalg_L$, is
defined by
$$
\begin{array}{lr}
\\{\amalg_K \vee^{\times} \amalg_L}(x)=\\\bigg\{<x,\alpha_K(x)\alpha_L(x), \beta_K(x)+\beta_L(x)-\beta_K(x)\beta_L(x)
>,f_{K \widetilde{\cup} L}(x)): x \in E \bigg\}
\end{array}
$$
where $f_{K \widetilde{\cup} L}(x)= f_K(x)\cup f_L(x)$.
\end{defn}
\begin{defn}\label{and-pro}
$\amalg_K, \amalg_L \in IFPS(U) $. Then AND-product $\amalg_K$ and
$\amalg_L$, denoted by $ \amalg_K\wedge^{\times} \amalg_L$,  is
defined by
$$
\begin{array}{l}
\\{\amalg_K\wedge^{\times} \amalg_L}(x)=\\\bigg\{<x,\alpha_K(x)\alpha_L(x), \beta_K(x)+\beta_L(x)-\beta_K(x)\beta_L(x)
>,f_{K \widetilde{\cap} L}(x)):x \in E \bigg\}
\end{array}
$$
where $f_{K \widetilde{\cap} L}(x)= f_K(x)\cap f_L(x)$.
\end{defn}
\begin{prop}
Let $\amalg_K,\amalg_L,\amalg_M\in IFPS(U) $. Then
\begin{enumerate}[i.]
\item $ \amalg_K \wedge^{\times}  \amalg_\Phi = \amalg_\Phi$
\item $ \amalg_K \wedge^{\times}  \amalg_{\tilde E} = \amalg_K  $
\item $ \amalg_K \wedge^{\times}  \amalg_L= \amalg_L \wedge^{\times}  \amalg_K $
\item $ \amalg_K \vee^{\times}  \amalg_L= \amalg_L \vee^{\times}  \amalg_K $
\item $(\amalg_K \wedge^{\times}  \amalg_L) \wedge^{\times}  \amalg_M= \amalg_K \wedge^{\times}  (\amalg_L\wedge^{\times}  \amalg_M) $
\item $(\amalg_K \vee^{\times}  \amalg_L) \vee^{\times}  \amalg_M= \amalg_K \vee^{\times}  (\amalg_L\vee^{\times}  \amalg_M) $
\end{enumerate}
\end{prop}



\section{Intuitionistic FP-soft decision
making method}

In this section, we have defined a reduced intuitionistic fuzzy set
of an intuitionistic FP-soft set, that produce an intuitionistic
fuzzy set from an intuitionistic FP-soft set. We then have defined a
reduced fuzzy set of an intuitionistic fuzzy set, that produce a
fuzzy set from an intuitionistic fuzzy set. These sets present an
adjustable approach to intuitionistic FP-soft sets based decision
making problems.
\begin{defn}
Let $\amalg_K$ be an intuitionistic FP-soft set. Then, a reduced
intuitionistic fuzzy set of  $\amalg_K$, denoted by $K_{rif}$,
defined as follows
$$
K_{rif}=\bigg\{<u,\alpha_{K_{rif}}(u),\beta_{K_{rif}}(u)> : u\in
U\bigg\}
$$
where
$$
\alpha_{K_{rif}}:U \rightarrow [0,1],\,\,\,
\alpha_{K_{rif}}(u)=\frac{1}{|U|}\sum_{x\in E,u\in
U}\alpha_K(x)\chi_{f_K(x)}(u)
$$
$$
\beta_{K_{rif}}:U \rightarrow [0,1],\,\,\,
\beta_{K_{rif}}(u)=\frac{1}{|U|}\sum_{x\in E,u\in
U}\beta_K(x)\chi_{f_K(x)}(u)
$$
where
$$
\chi_{f_K(x)}(u)=\left\{
\begin{array}{ll}
1, & u\in f_K(x)\\
0, & u \notin f_K(x)
\end{array}\right.
$$
Here, $\alpha_{K_{rif}}$ and  $\beta_{K_{rif}}$ are called
reduced-set operators of $K_{rif}$. It is clear that $K_{rif}$  is
an intuitionistic fuzzy set over $U$.
\end{defn}
\begin{defn}
$\amalg_K\in IFPS(U)$ and $K_{rif}$ be reduced intuitionistic fuzzy
set of $\amalg_K$. Then, a reduced fuzzy set of $K_{rif}$ is a fuzzy
set over $U$, denoted by $K_{rf}$, defined as follows
$$
K_{rf}=\bigg\{\mu_{K_{rf}}(u)/u:u\in U\bigg\}
$$
{ where}
$$
\mu_{K_{rf}}: U\rightarrow [0,1],\,\,\,
\mu_{K_{rf}}(u)=\alpha_{K_{rif}}(u)(1-\beta_{K_{rif}}(u))
$$
\end{defn}
Now, we construct an intuitionistic FP-soft decision making method
by the following algorithm to produce a decision fuzzy set from a
crisp set of the alternatives.

According to the problem, decision maker;
\begin{enumerate}[i.]
\item  constructs a feasible intuitionistic fuzzy subsets ${K}$  over
the parameters set $E$,
\item  constructs an intuitionistic FP-soft set $\amalg_K$  over the alternatives set $U$,
\item  computes the reduced intuitionistic fuzzy
set $K_{rif}$ of $\amalg_K$,
\item computes the reduced fuzzy set $K_{rf}$ of $K_{rif}$,
\item chooses the element of $K_{rf}$ that has maximum membership
degree.
\end{enumerate}
Now, we can give an example for the intuitionistic FP-soft decision
making method. Some of it is quoted from example in \cite{cag-09a}.
\begin{exmp}
Assume that a company wants to fill a position. There are 5
candidates who fill in a form in order to apply formally for the
position. There is a decision maker (DM), that is from the
department of human resources. He want to interview the candidates,
but it is very difficult to make it all of them. Therefore, by using
the intuitionistic FP-soft decision making method, the number of
candidates are reduced to a suitable one. Assume that the set of
candidates $U=\{u_1, u_2, u_3, u_4, u_5\}$  which may be
characterized by a set of parameters  $E=\{a_1, a_2, a_3, a_4\}$.
For $i=1,2,3,4$ the parameters $a_i$ stand for ''experience",
''computer knowledge", ''training" and ''young age", respectively.
Now, we can apply the method as follows:
\begin{enumerate}[\emph{Step} i.]
\item
Assume that DM constructs a feasible intuitionistic fuzzy subsets
${K}$ over the parameters set $E$ as follows;
$$
K= \{<x_1,0.7,0.3>, <x_2,0.2,0.5 >,<x_3,0.5,0.5 >, <x_4,0.6,0.3>\}
$$
\item
DM constructs an intuitionistic FP-soft set $\amalg_K$  over the
alternatives set $U$  as follows;
$$
\begin{array}{rl}
\amalg_K= & \bigg\{(<x_1,0.7,0.3 >,\{u_1, u_2,  u_4\}),
(<x_2,0.2,0.5>, U), (<x_3,\\&0.5,0.5 >, \{u_1,
u_2,u_4\}),(<x_4,0.6,0.3 >,\{ u_2, u_3\})\bigg\} \\
\end{array}
$$
\item
DM computes the reduced intuitionistic fuzzy set $K_{rif}$ of
$\amalg_K$ as follows;
$$
\begin{array}{rl}
K_{rif}= &\bigg\{(<u_1,0.28,0.26 >, <u_2,0.40,0.32>, <u_3,0.16,0.16
>,\\&<u_4,0.28,0.32 >, <u_5,0.04,0.10 >\bigg\}
\end{array}
$$

\item
DM computes the reduced fuzzy set $K_{rf}$ of $K_{rif}$ as follows;
$$
\begin{array}{rl}
K_{rf} =& \bigg\{0.2072/u_1, 0.2720/u_2, 0.1344/u_3, 0.1904/u_4,
0.0360/u_5\bigg\} \\
\end{array}
$$
\item
Finally, DM chooses  $u_2$ for the position from $K_{rf}$ since it
has the maximum degree 0.2720 among the others.
\end{enumerate}
\end{exmp}
Note that this decision making method can be applied for group
decision making easily with help of the Definition  \ref{or-sum},
Definition \ref{and-sum} Definition \ref{or-pro} and Definition
\ref{and-pro}.


\section{Conclusion}

In this paper, we first defined intuitionistic FP-soft sets and
their operations. We then presented the decision making method on
the intuitionistic FP-soft set theory. Finally, we provided an
example that demonstrated that the decision making method can
successfully work. It can be applied to problems of many fields that
contain uncertainty such as computer science, game theory, and so
on.


\end{document}